\documentclass[12pt]{amsart}

\usepackage{amsmath,amssymb,epic,lscape}

\newtheorem{theorem}{Theorem}[section]
\newtheorem{proposition}[theorem]{Proposition}

\theoremstyle{definition}

\theoremstyle{remark}

\numberwithin{equation}{section}

\def\DJ{{\hbox{D\kern-.8em\raise.15ex\hbox{--}\kern.35em}}}
\def\DJo{$\;$\kern-.4em
    \hbox{D\kern-.8em\raise.15ex\hbox{--}\kern.35em okovi\'c}}

\def\al{{\alpha}}

\def\vf{{\varphi}}
\def\la{{\lambda}}

\def\bZ{{\mbox{\bf Z}}}

\renewcommand{\subjclassname}{\textup{2000} Mathematics Subject
Classification }

\begin{document}

\title[Skew-Hadamard matrices]
{Skew-Hadamard matrices of orders $436$, $580$ and $988$ exist}

\author[D.\v{Z}. \DJ okovi\'{c}]
{Dragomir \v{Z}. \DJ okovi\'{c}}

\address{Department of Pure Mathematics, University of Waterloo,
Waterloo, Ontario, N2L 3G1, Canada}

\email{djokovic@uwaterloo.ca}

\thanks{
The author was supported by the NSERC Grant A-5285.}


\date{}

\begin{abstract}
We construct two difference families on each of the cyclic groups
of order $109$, $145$ and $247$, and use them to construct
skew-Hadamard matrices of orders $436$, $580$ and $988$.
Such difference families and matrices are constructed here
for the first time. The matrices are constructed by using
the Goethals--Seidel array.
\end{abstract}

\maketitle
\subjclassname{ 05B20, 05B30 }
\vskip5mm

\section{Introduction}

Recall that a Hadamard matrix of order $m$ is a $\{\pm1\}$-matrix
$A$ of size $m\times m$ such that $AA^T=mI_m$, where
$T$ denotes the transpose and $I_m$ the identity matrix.
A skew-Hadamard matrix is a Hadamard matrix $A$ such that
$A-I_m$ is a skew-symmetric matrix.

In this note we construct two non-equivalent skew-Hadamard
matrices for each of the orders $4\cdot109=436$, $4\cdot145=580$
and $4\cdot247=988$. They are constructed from cyclic difference 
families having suitable parameter sets and having another
special feature, namely one of the blocks is of skew type.
(The exact meaning of this term is explained below.)
We refer the reader to \cite{DZ} for the list of orders of
skew-Hadamard matrices constructed in our previous papers.
Our method of construction of concrete difference families is
described in those papers and will not be elaborated upon here.

The list of odd integers $n<300$ for which no skew-Hadamard matrix
of order $4n$ is known is now reduced to the following $39$ integers:
\begin{eqnarray*}
&& 69,89,101,107,119,149,153,167,177,179,191,193,201, \\
&& 205,209,213,223,225,229,233,235,239,245,249,251,253, \\
&& 257,259,261,265,269,275,277,283,285,287,289,295,299.
\end{eqnarray*}

We denote by $\bZ_n$ the ring, as well as the additive group,
of integers modulo $n$. By $\bZ_n^*$ we denote the group of invertible
elements of the ring $\bZ_n$. Its order is equal to $\vf(n)$
where $\vf$ is the Euler function.
We racall that a subset $X\subseteq\bZ_n$ is of \emph{skew type} 
if $X$ and $-X$ are disjoint and their union is $\bZ_n\setminus\{0\}$.

We construct our Hadamard matrices by using the Goethals--Seidel array
\[
\left[ \begin{array}{cccc}
		P_1 & P_2R & P_3R & P_4R \\
		-P_2R & P_1 & -P_4^TR & P_3^TR \\
		-P_3R & P_4^TR & P_1 & -P_2^TR \\
		-P_4R & -P_3^TR & P_2^TR & P_1
\end{array} \right] \]
consisting of 16 blocks of size $n$ by $n$.
As usual, $R$ denotes the matrix having ones on the back-diagonal
and all other entries zero. The matrices $P_i$ are usually circulant
$\{\pm1\}$-matrices constructed from difference families in
$\bZ_n$ having suitable parameters.

In the literature on Hadamard 
matrices it is customary to refer to difference families
as supplementary difference sets (SDS) and to employ
more elaborate and more informative notation by listing the
order $v$ of the underlying abelian group, the number of
sets in the family as well as their cardinals, and also the
parameter $\la$. We shall follow this practise.

\section{The case $n=109$}

In this section we set $n=109$. 
Observe that $n$ is a prime number and that $n-1=3\cdot36$.
Hence the multiplicative group $\bZ_n^*$
is cyclic and has the unique subgroup, namely $H=\{1,45,63\}$,
of order $3$. We shall use the same enumeration of the 36 
cosets $\al_i$, $0\le i\le 35$, of $H$ in $\bZ_n^*$ as in our 
computer program. Thus we impose the condition that
$\al_{2i+1}=-\al_{2i} \pmod{109}$ for $0\le i\le 17$.
For even indices we have
\[
\begin{array}{lllll}
\al_0=H, \quad & \al_2=2H, \quad & \al_4=3H, \quad & \al_6=4H, \quad & \al_8=5H, \\
\al_{10}=6H, & \al_{12}=8H, & \al_{14}=9H, & \al_{16}=10H, & \al_{18}=11H, \\
\al_{20}=13H, & \al_{22}=15H, & \al_{24}=16H, & \al_{26}=18H, & \al_{28}=20H, \\
\al_{30}=23H, & \al_{32}=25H, & \al_{34}=30H. &&
\end{array}
\]

Define the eight subsets of $\bZ_n$ by
\[ A_k = \bigcup_{i\in J_k} \al_i,\quad 
B_k = \bigcup_{i\in L_k} \al_i, \]
where
\begin{eqnarray*}
J_1 &=& \{ 0,2,5,7,8,10,12,15,16,19,20,23,24,26,29,30,33,34 \}, \\
J_2 &=& \{ 4,5,6,7,11,15,18,19,20,22,25,30,32,33,35 \}, \\
J_3 &=& \{ 0,1,5,6,9,10,11,14,17,20,24,26,27,28,29,31,32 \}, \\
J_4 &=& \{ 0,3,4,6,7,9,10,12,13,22,24,25,26,27,28,29,31,33,35 \}, \\
L_1 &=& \{ 0,3,5,6,9,10,13,15,16,18,21,23,24,27,28,30,32,34 \}, \\
L_2 &=& \{ 1,4,5,6,7,8,10,11,12,13,16,19,20,21,25,27,30,31,34 \}, \\
L_3 &=& \{ 0,1,2,3,4,6,10,12,13,15,16,17,20,22,23,24,25,29, \\
&& \quad 34,35 \}, \\
L_4 &=& \{ 0,1,2,3,5,8,9,12,13,14,15,16,19,21,23,25,26,28,29, \\
&& \quad 33,34 \}.
\end{eqnarray*}
Explicitly, we have
\begin{eqnarray*}
A_1 &=& \{ 1,2,5,6,7,8,10,13,14,16,17,18,23,27,29,30,31,32,33, \\
&& \quad 36,37,38,40,42,44,45,47,48,50,51,52,54,56,60,63, \\
&& \quad 66,68,70,74,75,81,83,84,85,87,88,89,90,94,97,98,100, \\
&& \quad 105,106 \},
\end{eqnarray*}
\begin{eqnarray*}
A_2 &=& \{ 3,4,11,13,15,21,23,25,26,29,31,32,34,35,38,39,40,43, \\
&& \quad 49,50,54,56,57,58,59,60,67,70,71,72,73,74,75,79,80, \\
&& \quad 82,83,84,87,93,98,100,103,105,106 \},
\end{eqnarray*}
\begin{eqnarray*}
A_3 &=& \{ 1,4,6,9,12,13,16,18,20,22,24,25,27,28,29,34,35,40,44, \\
&& \quad 45,46,47,48,49,51,52,55,56,57,58,61,62,63,64,65,66, \\
&& \quad 71,77,78,81,83,86,89,91,95,99,102,103,104,106,108 \},
\end{eqnarray*}
\begin{eqnarray*}
A_4 &=& \{ 1,3,4,6,8,12,15,16,18,19,20,21,26,27,28,33,34,38,41, \\
&& \quad 43,44,45,47,48,51,52,55,60,61,62,63,65,66,67,68,71, \\
&& \quad 72,73,74,75,76,77,79,80,81,82,84,86,89,91,92,93, \\
&& \quad 101,102,104,105,107 \}
\end{eqnarray*}
and
\begin{eqnarray*}
B_1 &=& \{ 1,4,6,10,11,12,14,16,19,20,23,25,27,28,29,30,31,32,34, \\
&& \quad 35,36,37,39,41,42,45,49,51,52,53,54,59,61,62,63,65,66, \\
&& \quad 69,71,76,83,85,87,88,91,92,94,96,100,101,102,104, \\
&& \quad 106,107 \},
\end{eqnarray*}
\begin{eqnarray*}
B_2 &=& \{ 3,4,5,6,7,8,10,13,14,23,26,29,30,32,33,34,37,38,40,41, \\
&& \quad 42,43,46,50,51,52,53,54,55,56,57,58,62,64,65,68,69, \\
&& \quad 70,71,75,76,77,80,82,83,85,86,91,93,96,97,98,101, \\
&& \quad 103,105,106,108 \},
\end{eqnarray*}
\begin{eqnarray*}
B_3 &=& \{ 1,2,3,4,6,8,10,13,14,15,16,17,19,21,24,26,27,30,31,33, \\
&& \quad 34,36,37,40,41,42,43,45,46,48,51,52,56,63,64,66,67,68, \\
&& \quad 71,72,73,76,79,80,81,82,85,87,88,89,90,92,93,94,95,99, \\
&& \quad 100,101,107,108 \},
\end{eqnarray*}
\begin{eqnarray*}
B_4 &=& \{ 1,2,5,7,8,9,10,12,14,17,18,19,20,22,28,29,30,31,33,36, \\
&& \quad 37,41,42,43,44,45,46,47,48,50,53,60,61,63,64,68,69,70, \\
&& \quad 74,76,78,81,82,83,84,85,87,88,89,90,92,93,94,96,97,98, \\
&& \quad 100,101,102,104,106,107,108 \}.
\end{eqnarray*}

Our heuristic search procedure was greatly simplified by searching
only through the special kind of subsets, namely those that are
union of cosets of $H$.

The cardinals of the above subsets are:
\begin{eqnarray*}
&& |A_1|=54,\ |A_2|=45,\ |A_3|=51,\ |A_4|=57, \\
&& |B_1|=54,\ |B_2|=57,\ |B_3|=60,\ |B_4|=63.
\end{eqnarray*}

By using a computer, one can easily verify the following
\begin{proposition} \label{stav1}
The quadruples $(A_1,A_2,A_3,A_4)$ and $(B_1,B_2,B_3,B_4)$ are
$4-(109;54,45,51,57;98)$ and $4-(109;54,57,60,63;125)$
supplementary difference sets in $\bZ_{109}$, respectively.
\end{proposition}

Any SDS is associated with a decomposition of $4n$ into sum of squares.
In the above two cases these decompositions are:
\[ 4n = \sum_{i=1}^4 (n-2|A_i|)^2 = \sum_{i=1}^4 (n-2|B_i|)^2,  \]
i.e., $436=1^2+19^2+7^2+5^2$ and $436=1^2+5^2+11^2+17^2$, respectively.

For any subset $X\subseteq\bZ_n$ let
\[ a_X=(a_0,a_1,\ldots,a_{n-1}) \]
be the $\{\pm1\}$-row vector such that $a_i=-1$ iff $i\in X$.
We denote by $[X]$ the $n\times n$ circulant matrix having
$a_X$ as its first row.

We can now substitute the matrices $[A_i]$ or $[B_i]$ for
$P_i$ in the Goethals--Seidel template
to obtain two Hadamard matrices of order $4n$.
		
\begin{proposition} \label{stav2}
The two Hadamard matrices of order $436$ constructed above from the
two SDS's of Proposition \ref{stav1} are skew-Hadamard matrices.
\end{proposition}
\begin{proof}
This follows from the fact that the sets $A_1$ and $B_1$ are of skew type.
\end{proof}

\section{The case $n=145$}

In this section we set $n=145$. The order of
$\bZ_n^*$ is $\vf(n)=\vf(5\cdot29)=4\cdot28$. Hence $\bZ_n^*$ has a
unique subgroup $H$ of order $7$, namely
\[ H=\{ 1,16,36,81,111,136,141 \}. \]
It acts on this ring by multiplication. There are $25$ orbits.
We enumerate the $24$ nonzero orbits by setting:
\begin{eqnarray*}
\al_0 &=& H, \\
\al_2 &=& \{ 2,17,32,72,77,127,137 \}, \\
\al_4 &=& \{ 3,43,48,98,108,118,133 \}, \\
\al_6 &=& \{ 6,51,71,86,91,96,121 \}, \\
\al_8 &=& \{ 7,52,82,107,112,117,132 \}, \\
\al_{10} &=& \{ 11,21,31,46,61,101,106 \}, \\
\al_{12} &=& \{ 14,19,69,79,89,104,119 \}, \\
\al_{14} &=& \{ 22,42,57,62,67,92,122 \}, \\
\al_{16} &=& \{ 5,35,80,100,115,120,125 \}, \\
\al_{18} &=& \{ 10,15,55,70,85,95,105 \}, \\
\al_{20} &=& \{29\}, \\
\al_{22} &=& \{58\}
\end{eqnarray*}
and $\al_{2i+1}=-\al_{2i} \pmod{145}$ for $0\le i\le11$.

Define the eight subsets of $\bZ_n$ by
\[ A_k = \bigcup_{i\in J_k} \al_i,\quad 
B_k = \bigcup_{i\in L_k} \al_i, \]
where
\begin{eqnarray*}
J_1 &=& \{ 1,2,4,7,9,10,13,14,16,19,20,22 \}, \\
J_2 &=& \{ 0,2,4,7,10,11,14,18,19,20,21,22 \}, \\
J_3 &=& \{ 1,3,6,9,12,13,14,17,19,20,21,22,23 \}, \\
J_4 &=& \{ 2,3,5,6,7,9,12,13,15,16,19,20,21,22,23 \}, \\
L_1 &=& \{ 1,3,5,7,8,10,13,15,17,18,20,23 \}, \\
L_2 &=& \{ 0,1,2,4,12,13,15,18,19,20,22,23 \}, \\
L_3 &=& \{ 4,5,6,7,9,12,14,16,18,20,21,22,23 \}, \\
L_4 &=& \{ 2,6,9,10,11,12,13,14,15,17,18,20,21,22,23 \}.
\end{eqnarray*}
Explicitly we have
\begin{eqnarray*}
A_1 &=& \{ 2,3,4,5,9,11,13,17,21,22,24,26,28,29,31,32,33,34,35, \\
&& \quad 38,40,41,42,43,46,48,49,50,54,56,57,58,59,60,61,62, \\
&& \quad 63,64,66,67,72,74,75,76,77,80,90,92,93,94,98,100,101, \\
&& \quad 106,108,109,115,118,120,122,125,126,127,129,130,131, \\
&& \quad 133,135,137,138,139,144  \},
\end{eqnarray*}
\begin{eqnarray*}
A_2 &=& \{ 1,2,3,10,11,15,16,17,21,22,24,29,31,32,36,39,40,42,43, \\
&& \quad 44,46,48,49,50,54,55,57,58,59,60,61,62,67,70,72,74, \\
&& \quad 75,77,81,84,85,90,92,94,95,98,99,101,105,106,108,111, \\
&& \quad 114,116,118,122,124,127,130,133,134,135,136,137, \\
&& \quad 139,141 \},
\end{eqnarray*}
\begin{eqnarray*}
A_3 &=& \{ 4,6,8,9,13,14,18,19,20,22,25,26,28,29,30,33,34,38,40, \\
&& \quad 41,42,45,50,51,56,57,58,60,62,63,64,65,66,67,68,69,71, \\
&& \quad 73,75,76,79,86,87,89,90,91,92,93,96,104,109,110,113, \\
&& \quad 116,119,121,122,126,128,129,130,131,135,138,140, \\
&& \quad 143,144 \},
\end{eqnarray*}
\begin{eqnarray*}
A_4 &=& \{ 2,5,6,8,12,13,14,17,18,19,23,24,26,27,28,29,32,33,35, \\
&& \quad 37,38,40,41,47,49,50,51,53,54,56,58,59,60,63,66,68,69, \\
&& \quad 71,72,73,74,75,76,77,78,79,80,83,86,87,88,89,90,91,93, \\
&& \quad 94,96,97,100,102,103,104,113,115,116,119,120,121,123, \\
&& \quad 125,126,127,128,130,131,135,137,138,139,142,143 \}
\end{eqnarray*}
and
\begin{eqnarray*}
B_1 &=& \{ 4,7,8,9,10,11,12,15,18,20,21,23,24,25,26,27,29,30,31, \\
&& \quad 34,37,41,45,46,47,49,52,53,54,55,56,59,61,64,65,66, \\
&& \quad 68,70,73,74,76,78,82,83,85,87,88,94,95,97,101,102, \\
&& \quad 103,105,106,107,109,110,112,113,117,123,126,128, \\
&& \quad 129,131,132,139,140,142,143,144 \}, \\
\end{eqnarray*}
\begin{eqnarray*}
B_2 &=& \{ 1,2,3,4,9,10,14,15,16,17,19,23,26,29,32,34,36,40,41,43, \\
&& \quad 48,50,53,55,56,58,60,64,66,69,70,72,75,76,77,78,79,81, \\
&& \quad 83,85,87,88,89,90,95,98,103,104,105,108,109,111,118, \\
&& \quad 119,123,126,127,129,130,131,133,135,136,137,141,144  \}, \\
\end{eqnarray*}
\begin{eqnarray*}
B_3 &=& \{ 3,5,6,10,12,13,14,15,19,22,24,27,28,29,33,35,37,38,42, \\
&& \quad 43,47,48,49,51,54,55,57,58,59,62,63,67,69,70,71,74, \\
&& \quad 79,80,85,86,87,89,91,92,93,94,95,96,97,98,100,102, \\
&& \quad 104,105,108,115,116,118,119,120,121,122,125,133,138, \\
&& \quad 139,142  \}, \\
\end{eqnarray*}
\begin{eqnarray*}
B_4 &=& \{ 2,6,10,11,13,14,15,17,19,20,21,22,23,25,26,28,29,30,31, \\
&& \quad 32,33,38,39,41,42,44,45,46,51,53,55,56,57,58,61,62,63, \\
&& \quad 65,66,67,69,70,71,72,76,77,78,79,83,84,85,86,87,88,89, \\
&& \quad 91,92,93,95,96,99,101,103,104,105,106,110,114,116, \\
&& \quad 119,121,122,123,124,126,127,131,134,137,138,140 \}.
\end{eqnarray*}
The cardinals of these subsets are:
\begin{eqnarray*}
&& |A_1|=|B_1|=72,\ |A_2|=|B_2|=66, \\
&& |A_3|=|B_3|=67,\ |A_4|=|B_4|=81.
\end{eqnarray*}

By using a computer, one can easily verify the following
\begin{proposition} \label{stav3}
The quadruples $(A_1,A_2,A_3,A_4)$ and $(B_1,B_2,B_3,B_4)$ are
$4-(145;72,66,67,81;141)$ supplementary difference sets in $\bZ_{145}$.
\end{proposition}

Again our heuristic search procedure used only
the subsets which are unions of nonzero $H$-orbits. 
These two SDS's are not equivalent.
Indeed, we have verified that the sets $A_2$ and $B_2$ are
not equivalent under translations and multiplications by
elements of $\bZ_n^*$. As these two SDS's have the same parameters,
they are associated with the same decomposition of $4n$ into 
sum of squares: $580=1^2+13^2+11^2+17^2$.

Now observe that the subsets $A_1$ and $B_1$ are of skew type.
This is easy to verify since the index sets $J_1$ and $L_1$
meet each pair $\{2i,2i+1\}$, $0\le i\le11$, in exactly one element.
Hence the following holds.
		
\begin{proposition} \label{stav4}
The two Hadamard matrices of order $580$, which one can construct from
the two SDS's of Proposition \ref{stav3} by using the Goethals--Seidel
array, are skew-Hadamard matrices.
\end{proposition}

\section{The case $n=247$}

In this section we set $n=247$. As $n=13\cdot19$, we have
$|\bZ_n^*|=12\cdot18$. We select one of its cyclic subgroups
of order 9, namely
\[ H=\{ 1,9,16,55,61,81,139,144,235 \}. \]
It acts on the ring $\bZ_n$ by multiplication. There are $31$ orbits.
We enumerate the $30$ nonzero orbits by setting:
\begin{eqnarray*}
\al_0 &=& H, \\
\al_2 &=& \{ 2,18,31,32,41,110,122,162,233 \}, \\
\al_4 &=& \{ 3,27,48,165,170,183,185,211,243 \}, \\
\al_6 &=& \{ 5,28,45,58,80,158,187,201,226 \}, \\
\al_8 &=& \{ 6,54,83,93,96,119,123,175,239 \}, \\
\al_{10} &=& \{ 7,20,63,73,112,138,163,180,232 \}, \\
\al_{12} &=& \{ 10,56,69,90,116,127,155,160,205 \}, \\
\al_{14} &=& \{ 11,47,99,102,111,115,150,176,177 \}, \\
\al_{16} &=& \{ 13,52,65,78,91,117,143,208,221 \}, \\
\al_{18} &=& \{ 14,29,40,79,113,126,146,217,224 \}, \\
\al_{20} &=& \{ 17,25,43,49,140,142,153,194,225 \}, \\
\al_{22} &=& \{ 19,57,171 \}, \\
\al_{24} &=& \{ 33,34,37,50,59,86,98,141,203 \}, \\
\al_{26} &=& \{ 35,66,68,74,100,118,159,172,196 \}, \\
\al_{28} &=& \{ 38,95,114 \}
\end{eqnarray*}
and $\al_{2i+1}=-\al_{2i} \pmod{145}$ for $0\le i\le11$.

Define the eight subsets of $\bZ_n$ by
\[ A_k = \bigcup_{i\in J_k} \al_i,\quad 
B_k = \bigcup_{i\in L_k} \al_i, \]
where
\begin{eqnarray*}
J_1 &=& \{ 0,2,4,7,8,10,12,15,16,18,20,23,25,27,29 \}, \\
J_2 &=& \{ 0,2,7,9,11,12,14,15,16,18,20,22,26 \}, \\
J_3 &=& \{ 2,3,4,12,13,14,15,16,18,20,23,24,26,27,29 \}, \\
J_4 &=& \{ 0,3,4,6,10,11,12,14,18,19,20,22,25,29 \}, \\
L_1 &=& \{ 0,3,5,7,8,11,12,15,17,18,21,22,24,27,29 \}, \\
L_2 &=& \{ 3,5,6,8,11,13,14,15,16,19,26,27,29 \}, \\
L_3 &=& \{ 0,1,2,4,5,7,11,13,14,15,22,23,24,26,27 \}, \\
L_4 &=& \{ 0,3,8,9,10,11,13,17,19,24,25,27,28,29 \}.
\end{eqnarray*}
Explicitly, we have
\begin{eqnarray*}
A_1 &=& \{ 1,2,3,6,7,9,10,13,14,16,17,18,20,21,25,27,29,31,32, \\
&& \quad 40,41,43,44,46,48,49,51,52,54,55,56,60,61,63,65,69,70, \\
&& \quad 71,73,75,76,78,79,81,83,88,89,90,91,93,96,97,106,110, \\
&& \quad 112,113,116,117,119,122,123,126,127,129,132,133,136, \\
&& \quad 138,139,140,142,143,144,145,146,147,148,149,152,153, \\
&& \quad 155,160,161,162,163,165,167,170,173,175,179,180,181, \\
&& \quad 183,185,188,189,190,194,197,200,202,205,208,209,210, \\
&& \quad 211,212,213,214,217,219,221,223,224,225,228,232,235, \\
&& \quad 236,239,242,243 \},
\end{eqnarray*}
\begin{eqnarray*}
A_2 &=& \{ 1,2,8,9,10,11,13,14,15,16,17,18,19,21,25,29,31,32,35, \\
&& \quad 40,41,43,46,47,49,52,55,56,57,60,61,65,66,67,68,69,70, \\
&& \quad 71,72,74,78,79,81,84,89,90,91,97,99,100,102,109,110, \\
&& \quad 111,113,115,116,117,118,122,124,126,127,128,132,135, \\
&& \quad 136,139,140,142,143,144,145,146,148,150,151,153,154, \\
&& \quad 155,159,160,162,164,167,171,172,174,176,177,184,189, \\
&& \quad 193,194,196,200,202,205,208,217,219,221,223,224,225, \\
&& \quad 227,235,236,240,241,242 \},
\end{eqnarray*}
\begin{eqnarray*}
A_3 &=& \{ 2,3,10,11,13,14,17,18,24,25,27,29,31,32,33,34,35,37, \\
&& \quad 40,41,42,43,47,48,49,50,51,52,56,59,65,66,68,69,70,71, \\
&& \quad 74,75,76,78,79,85,86,87,88,90,91,92,97,98,99,100,102, \\
&& \quad 110,111,113,115,116,117,118,120,122,125,126,127,129, \\
&& \quad 131,132,133,136,137,140,141,142,143,145,146,147,148, \\
&& \quad 150,152,153,155,157,159,160,162,165,170,172,173,176, \\
&& \quad 177,178,179,181,183,185,190,191,194,196,200,203,205, \\
&& \quad 206,208,209,211,212,215,216,217,221,223,224,225,228, \\
&& \quad 229,236,237,243,245 \},
\end{eqnarray*}
\begin{eqnarray*}
A_4 &=& \{ 1,3,5,7,9,10,11,14,15,16,17,19,20,23,24,25,27,28,29, \\
&& \quad 30,40,43,44,45,47,48,49,55,56,57,58,61,63,67,69,73,79, \\
&& \quad 80,81,84,85,90,99,101,102,106,109,111,112,113,115,116, \\
&& \quad 121,125,126,127,133,134,135,137,138,139,140,142,144, \\
&& \quad 146,149,150,152,153,155,158,160,161,163,165,168,170, \\
&& \quad 171,174,176,177,180,183,184,185,187,188,194,197,201, \\
&& \quad 205,206,207,209,210,211,213,214,215,216,217,218,224, \\
&& \quad 225,226,227,229,232,233,235,240,243,245 \}
\end{eqnarray*}
and
\begin{eqnarray*}
B_1 &=& \{ 1,4,6,9,10,14,15,16,19,21,22,24,26,29,33,34,36,37,39, \\
&& \quad 40,46,50,51,53,54,55,56,57,59,60,61,62,64,67,69,70,71, \\
&& \quad 75,77,79,81,82,83,84,85,86,88,89,90,93,94,96,97,98, \\
&& \quad 104,105,107,109,113,116,119,123,125,126,127,129,130, \\
&& \quad 132,133,135,136,137,139,141,144,145,146,147,148,152, \\
&& \quad 155,156,160,167,169,171,173,174,175,179,181,182,184, \\
&& \quad 189,195,198,199,200,202,203,204,205,206,209,212,215, \\
&& \quad 216,217,219,220,222,224,227,229,230,234,235,236,239, \\
&& \quad 240,242,244,245 \},
\end{eqnarray*}
\begin{eqnarray*}
B_2 &=& \{ 4,5,6,11,13,15,23,24,28,30,35,36,42,45,47,51,52,54, \\
&& \quad 58,62,64,65,66,67,68,70,71,74,75,77,78,80,82,83,84, \\
&& \quad 85,87,88,91,92,93,96,97,99,100,101,102,109,111,115, \\
&& \quad 117,118,119,120,121,123,125,129,131,132,133,134,135, \\
&& \quad 136,137,143,145,147,148,150,152,157,158,159,168,172, \\
&& \quad 173,174,175,176,177,178,179,181,184,187,191,196,199, \\
&& \quad 200,201,206,207,208,209,212,215,216,218,220,221,226, \\
&& \quad 227,229,233,236,237,239,240,244,245 \},
\end{eqnarray*}
\begin{eqnarray*}
B_3 &=& \{ 1,2,3,4,9,11,12,15,16,18,19,21,27,31,32,33,34,35,36, \\
&& \quad 37,41,42,46,47,48,50,51,55,57,59,60,61,62,64,66,67, \\
&& \quad 68,70,71,74,75,76,77,81,82,84,86,87,88,89,92,97,98, \\
&& \quad 99,100,102,103,108,109,110,111,115,118,120,122,129, \\
&& \quad 131,132,135,136,139,141,144,145,147,148,150,157,159, \\
&& \quad 162,165,166,167,170,171,172,173,174,176,177,178,179, \\
&& \quad 181,183,184,185,186,189,190,191,192,196,199,200,202, \\
&& \quad 203,211,212,219,220,223,227,228,231,235,236,237,238, \\
&& \quad 240,242,243,244,246 \},
\end{eqnarray*}
\begin{eqnarray*}
B_4 &=& \{ 1,6,7,8,9,15,16,20,23,24,26,30,33,34,37,38,39,42,44, \\
&& \quad 50,51,54,55,59,61,63,67,72,73,75,81,83,84,85,86,87, \\
&& \quad 88,92,93,95,96,98,101,104,106,109,112,114,119,120,121, \\
&& \quad 123,124,125,128,129,130,131,133,134,135,137,138,139, \\
&& \quad 141,144,147,149,151,152,154,156,157,161,163,164,168, \\
&& \quad 169,173,174,175,178,179,180,181,182,184,188,191,193, \\
&& \quad 195,197,203,206,207,209,210,212,213,214,215,216,218, \\
&& \quad 227,229,232,233,234,235,237,239,240,241,245 \}.
\end{eqnarray*}
The cardinals of these subsets are:
\begin{eqnarray*}
&& |A_1|=|B_1|=123,\ |A_2|=|B_2|=111, \\
&& |A_3|=|B_3|=123,\ |A_4|=|B_4|=114.
\end{eqnarray*}

By using a computer, one can easily verify the following
\begin{proposition} \label{stav5}
The quadruples $(A_1,A_2,A_3,A_4)$ and $(B_1,B_2,B_3,B_4)$ are
$4-(247;123,111,123,114;224)$ supplementary difference sets in $\bZ_{247}$.
\end{proposition}

Again our heuristic search procedure used only
the subsets which are unions of nonzero $H$-orbits. 
The above two SDS's are not equivalent.
Indeed, we have verified that the sets $A_2$ and $B_2$ are
not equivalent under translations and multiplications by
elements of $\bZ_n^*$. As these two SDS's have the same parameters,
they are associated with the same decomposition of $4n$ into 
sum of squares: $988=1^2+1^2+19^2+25^2$.

Now observe that the subsets $A_1$ and $B_1$ are of skew type.
This is easy to verify since the index sets $J_1$ and $L_1$
meet each pair $\{2i,2i+1\}$, $0\le i\le15$, in exactly one element.
Hence the following holds.
		
\begin{proposition} \label{stav6}
The two Hadamard matrices of order $988$ constructed above from the
two SDS's of Proposition \ref{stav5} are skew-Hadamard matrices.
\end{proposition}

\end{document}